\documentclass[11pt]{article}
\usepackage{amsfonts}
\usepackage{a4wide,latexsym,amssymb,amsfonts,amsmath}
\usepackage{txfonts,color}
\newtheorem{theorem}{Theorem}[section]

\newtheorem{lemma}[theorem]{Lemma}
\newtheorem{corollary}[theorem]{Corollary}
\newtheorem{example}{Example}[section]
\newtheorem{construction}[theorem]{Construction}

\def\whitebox{{\hbox{\hskip 1pt
\vrule height 6pt depth 1.5pt \lower 1.5pt\vbox to 7.5pt{\hrule
width 3.2pt\vfill\hrule width 3.2pt} \vrule height 6pt depth 1.5pt
\hskip 1pt } }}
\def\qed{\ifhmode\allowbreak\else\nobreak\fi\hfill\quad\nobreak
\whitebox\medbreak}
\newcommand{\proof}{\noindent{\bf Proof.}\ }

\newcommand{\pf}{\noindent{\bf Proof.}\ }
\newcommand{\ignore}[1]{}

\newcommand{\A}{{\cal A}}

\newcommand{\Z}{{\mathbb Z}}

\begin{document}
\title{\bf Consecutive Detecting Arrays for Interaction Faults
\thanks{Correspondence to: Aiyuan Tao (taoaiyuan@126.com).
This work was supported by NSFC No. 11301342 and Natural Science Foundation of Shanghai No. 17ZR1419900.}}

\author{ Ce Shi$^1$, Ling Jiang$^2$ and Aiyuan  Tao$^3$ \\
\small$^1$  School of  Statistics and Mathematics \\
\small Shanghai Lixin University of Accounting and Finance, \
\small Shanghai 201209, China \\
\small$^2$ College  of Information Technology\\
\small Shanghai Ocean University, Shanghai 201306, China\\
\small$^3$  School of International Economics and Trade \\
\small Shanghai Lixin University of Accounting and Finance, \
\small Shanghai 201209, China \\
}

\maketitle

\begin{abstract}
The concept of detecting arrays was  developed to locate and detect interaction faults arising between the factors in
a component-based system during software testing. In this paper, we propose a family
of consecutive detecting arrays (CDAs) in which the interactions
between factors are considered to be ordered. CDAs can be used to generate test suites for
locating and detecting interaction faults between neighboring factors. We establish a general
criterion for measuring the optimality of  CDAs in terms of their size. Based on this optimality
criterion, the equivalence between optimum CDAs and  consecutive orthogonal arrays with
prescribed properties is explored. Using the advantages of this equivalence, a great number of
optimum CDAs are presented. In particular, the existence of optimum CDAs with few factors is  completely determined.

\vspace{0.2cm}
\noindent
 {\bf Keywords}: consecutive detecting arrays; consecutive covering arrays; consecutive orthogonal arrays;
 optimality; equivalence

\vspace{0.2cm}
\noindent
{\bf Mathematics Subject Classifications (2010)}: 05B15, 05B20, 62K15, 94C12
\end{abstract}

\section{Introduction}

Throughout this paper, let $I_n$ be the set of the first $n$ positive integers. An $N\times k$ array
with  entries from an alphabet $\{0,1,\cdots, v-1\}$ of size $v$ is said to be a covering array (CA)
(resp. orthogonal array (OA)) if each $N\times t$ subarray contains each $t$-tuple at least
(resp. exactly ) $\lambda$ times among its rows. It is denoted by CA$_\lambda(N;t,k,v)$ (resp. OA$_\lambda(N;t,k,v)$ ).
When $\lambda=1$, the notation CA$(N;t,k,v)$  (resp. OA$(t,k,v)$) is often used.

CAs are often  used to generate test suites in a component-based system, and have various applications
in software and hardware, circuit design, and so on. They present a useful alternative to exhaustive testing, and the
use of these arrays to construct test suites and their ability to dramatically reduce the testing burden are supported by
many empirical results \cite{KR2002,KW2004}. CAs have been studied extensively,
and a great number of methods and results have been reported \cite{CCK1999,CK2002, CMTW2006,Colbourn2004,Colbourn2008,HR2004,JM2018,JY2010,TMPS2017}.

If  the selection of the $t$ columns is restricted by
considering only consecutive columns, the ensuing problem is closely related to the CA framework. More specifically, Godbole {\em et al.} \cite{GKM2011} introduced the concept of a
consecutive covering array (CCA), in which the structure of the columns captures some linear progression of data (for example, data across a series of consecutive dates) or
data organized by consecutive proximity (for example, consecutive switches in a circuit). CCAs
(resp. consecutive orthogonal arrays (COAs)), denoted as CCA$(N;t,k,v)$ (resp. COA$_\lambda(t, k,v)$), are $N\times k$ arrays with
entries from a set $V$ of $v$ symbols such that each set of $t$ consecutive columns contains each $t$-tuple at least once
(resp. exactly $\lambda$ times) among its rows. The following array is a CCA$(9;2,21,3)$ over $\Z_3$.
$$
\left(
\begin{array}{ccccccccccccccccccccc}
  0 & 0 & 0  & 0 & 0 & 0 & 0 & 0 & 0 & 0  & 0 & 0 & 0 & 0 & 0 & 0 & 0 & 0 & 0 & 0& 0 \\
  0 & 1 & 0  & 1 & 0 & 1 & 0 & 1 & 0 & 1  & 0 & 1 & 0 & 1 & 0 & 1 & 0 & 1 & 0 & 1& 0 \\
  0 & 2 & 0  & 2 & 0 & 2 & 0 & 2 & 0 & 2  & 0 & 2 & 0 & 2 & 0 & 2 & 0 & 2 & 0 & 2& 0 \\
  1 & 0 & 1  & 0 & 1 & 0 & 1 & 0 & 1 & 0  & 1 & 0 & 1 & 0 & 1 & 0 & 1 & 0 & 1 & 0& 1 \\
  1 & 1 & 1  & 1 & 1 & 1 & 1 & 1 & 1 & 1  & 1 & 1 & 1 & 1 & 1 & 1 & 1 & 1 & 1 & 1& 1 \\
  1 & 2 & 1  & 2 & 1 & 2 & 1 & 2 & 1 & 2  & 1 & 2 & 1 & 2 & 1 & 2 & 1 & 2 & 1 & 2& 1 \\
  2 & 0 & 2  & 0 & 2 & 0 & 2 & 0 & 2 & 0  & 2 & 0 & 2 & 0 & 2 & 0 & 2 & 0 & 2 & 0& 2 \\
  2 & 1 & 2  & 1 & 2 & 1 & 2 & 1 & 2 & 1  & 2 & 1 & 2 & 1 & 2 & 1 & 2 & 1 & 2 & 1& 2 \\
  2 & 2 & 2  & 2 & 2 & 2 & 2 & 2 & 2 & 2  & 2 & 2 & 2 & 2 & 2 & 2 & 2 & 2 & 2 & 2& 2
\end{array}
\right).
$$

CCAs represent a similar family to CAs. The analogy between CAs and CCAs
is almost obvious. Just as CAs can be used for combinatorial testing, CCAs can be used to generate test suites
for combinatorial testing of neighboring factors, such as in circuit testing, signal processing, and so on \cite{WNXS2007}. Although CAs can also be used to test such software systems,  they are larger than CCAs for fixed $t,k,v$. For example, the best upper bound on
the size of CA$(N;2,21,3)$ is 16, but the above example gives us a CCA with only 9 rows. In \cite{GKM2011}, Godbole {\em et al.} focused on CCA$(N;t,k,2)$
using a probabilistic approach  based on an appropriate  Markov chain method. This allowed them to determine the probability distribution function of a random variable that enumerates the number of uncovered consecutive $t$-subarrays in the case of an $N\times k$ binary array obtained by realizing $kn$ Bernoulli
variables. The more general problem of CCA$(N;t,k,v)$, i.e., establishing the probability distribution function of
the random variable enumerating the uncovered consecutive $t$-subarrays with the Markov chain method, was considered in \cite{GKM2010}.
Optimal  consecutive covering arrays are studied and generalized to optimal cyclic consecutive covering arrays, and many results
are given in \cite{RMS2018}.

When using  CAs to generate test suites, the columns of the CA represent factors affecting the response and the entries within columns indicate the settings
or values for that factor. The rows then represent the tests to be run. Thus, testing with a CA
can indicate the presence or absence of interaction faults for up to $t$ factors.  This constitutes a
valuable step in the process of screening a system for interaction faults prior to its release.
However, the location and magnitude of the interactions causing the faults may be far from clear.
In practical terms,  tests that reveal the location of the interaction faults are of considerable interest.  For this,
Colbourn and McClary formalized the problem of the nonadaptive location of interaction faults under the
hypothesis that the system contains (at most) $d$ faults, each involving (at most) $t$ interacting factors \cite{CM2008}. They proposed the notion of detecting arrays  to solve this problem.

Let  $A=(a_{ij}) \ (i \in I_N, j\in I_k)$ be an $N\times k$ array with entries from an alphabet $V$ of size $v$. Each $t$-set of
columns with $t$-tuples of values for those columns is called a $t$-way interaction, denoted by $T=\{(j_r,x_r): x_r\in V, 1\leq r\leq t  \}$,
where  $1\le j_1<j_2 < \cdots < j_t \le k$. Write $\rho(A, T)$ for the set of indices of rows of $A$ that cover $T$, i.e.,  $\rho(A, T)= \{i: a_{ij_r}=x_r, 1\le r\le t\}$.
For an arbitrary set $\mathcal {T}$ of interactions, define $\rho (A,\mathcal {T})=\cup_{T\in {\cal T}}\rho (A,T)$. Furthermore, suppose that $A=(a_{ij})\ (i \in I_N, j\in I_k)$ is a CA$(N; t, k, v )$ over $V$. Write
$\mathcal {I}_t$ for the set of all $t$-way interactions of $A$. For any ${\cal T}\subseteq \mathcal {I}_t$  with $|{\cal T}| = d$
and any $T\in \mathcal {I}_t$, if we have
$$\rho(A, T)\subseteq \rho(A, {\cal T})
\Leftrightarrow T\in {\cal T},$$

\noindent then the array $A$  is called a $(d,t)$-{\em detecting array}
 (DA), denoted by $(d,t)$-DA$(N;k,v)$.

Similar to CAs,  CCAs can generate test suites for combinatorial testing of neighboring factors to indicate
the presence or absence of faulty interactions, but they cannot identify and determine these interactions from the outcome
of tests.  Although DAs can be used to locate and detect interaction faults between neighboring
factors,  they are not well adapted for this kind of software testing. For example, some optimum DAs do
not exist for fixed $t,k,v$ \cite{STY2012,SY2014}, but the optimum arrays for locating and detecting interaction
faults between neighboring factors may exist. Moreover, combinatorial testing of neighboring factors only considers
consecutive interactions, rather than arbitrary  interactions. In an attempt to solve this problem, we propose a similar family of DAs called consecutive detecting arrays (CDAs), in which the interactions between factors are considered to be ordered.

The remainder of this paper is organized as follows. The concept of CDAs is described in Section $2$. The necessity for the existence of  CDAs is also discussed in this section.
In Section $3$, we establish a general criterion for measuring the optimality of a CDA in terms of its size, and explore the
equivalence between optimum CDAs and  COAs with prescribed properties. Based on the equivalence
outlined in Section $3$, some  constructions and existence results are presented in Section $4$. Finally, Section $5$ summarizes our concluding remarks.

\section{Consecutive Detecting Arrays}

This section explains the notion of CDAs. The necessity for
the existence of CDAs  is also discussed. To aid this discussion, a consecutive $t$-way interaction is defined below.

A consecutive $t$-way interaction is denoted as $T=\{((i,x_i),(i+1,x_{i+1}),\cdots,(i+t-1,x_{i+t-1}))\}$,
where $1\leq i\leq k-t+1$, $x_r\in V$ for $r=i,i+1,\cdots, i+t-1$. Obviously, there are a total of $(k-t+1)v^t$ consecutive $t$-way
interactions for $k$ neighboring factors. To locate and detect  interaction faults between neighboring factors, it is only
necessary to identify the consecutive interaction faults from the outcomes of the tests. Thus, the notion of  CDAs comes from modifying the notion of DAs.

Suppose  that $A=(a_{ij})\ (i \in I_N, j\in I_k)$ is a CCA$(N; t, k, v )$ over $V$. Write $C\mathcal {I}_t$ for the set of all consecutive  $t$-way interactions of $A$. For any ${\cal T}\subseteq C \mathcal {I}_t$  with $|{\cal T}| = d$
and any $T\in C \mathcal {I}_t$, if we have
$$\rho(A, T)\subseteq \rho(A, {\cal T})
\Leftrightarrow T\in {\cal T},$$

\noindent then the array $A$  is called a $(d,t)$-{\em  consecutive  detecting array}, denoted by $(d,t)$-CDA$(N;k,v)$.

Clearly, a $(d,t)$-DA$(N;k,v)$ must be a $(d,t)$-CDA$(N;k,v)$, but the converse is not always true. It is straightforward that $T\in {\cal T}$ implies $\rho(A, T)\subseteq \rho(A, {\cal T})$.
Hence, the condition $\rho(A, T)\subseteq \rho(A, {\cal T})\Leftrightarrow T\in {\cal T},$ is satisfied
if $T\not \in {\cal T}\Rightarrow\rho(A, T)\not \subseteq \rho(A, {\cal T})$. We will make extensive
use of this simple fact in the following. As well as DAs, there are some
admissible parameters for the existence of CDAs. We restrict our discussion to nontrivial parameters.  As there are exactly $(k-t+1)v^t$
possible consecutive $t$-way interactions, we treat $(d,t)$-CDA only when $1\leq d\leq (k-t+1)v^t$. When $k<t$
and $d>0$, no $(d,t)$-CDAs can exist. If $k=t$, we can form an array consisting of all $t$-tuples. Hence, we only treat
cases with $k>t$. Finally, we require $v > 1$  to avoid factors that take on unique levels. The following results
 for DAs were given in \cite{CM2008}.

\begin{lemma}
 Suppose that $A$ is a $(d,t)$-DA$(N;k,v)$. Then
\begin{enumerate}
\item $d<v$.
\item $A$ is also a $(s,t)$-DA$(N;k,v)$, where $1\leq s\leq d-1$.
\end{enumerate}
\end{lemma}

As a direct consequence of DAs,  the following  lemmas can be easily obtained. We state them for later use.

\begin{lemma}\label{NC-CDA}
If a $(d,t)$-CDA$(N;k,v)$ exists, then $d<v$.
\end{lemma}

\begin{lemma}\label{(d,t)-DA-(d-1,t)-DA}
Let $A$ be a $(d,t)$-CDA$(N;k,v)$ that exists. Then, $A$ is also a $(s,t)$-CDA$(N;k,v)$, where $1\leq s\leq d-1$.
\end{lemma}

By definition, a $(d,t)$-CDA is actually a special CCA of strength $t$. The significance of using the CDA  to generate test suites
 is that any set of $d$   consecutive $t$-way interaction faults can be determined from the outcomes. Further, if there are more than $d$ consecutive $t$-way
interactions causing the faults, this can also be detected.  For details, see the application of DAs in \cite{CM2008}. As the rows of
a CDA stand for tests, the CDA of minimum size when other parameters are fixed is of considerable interest. The minimum $N$ for which a $(d,t)$-CDA$(N;k,v)$ exists is
referred to as the {\em consecutive detecting array number} (CDAN), denoted by $(d,t)$-CDAN$(k,v)$. A $(d,t)$-CDA$(N;k,v)$ with
$N=(d,t)$-\mbox{CDAN}$(k,v)$ is said to
be optimum.  In the next section, we derive a lower bound for the function $(d,t)$-CDAN$(k,v)$.

\section{Optimality Criterion and Combinatorial Description}

The objective of this section is to establish a lower bound on the size of $(d,t)$-CDA$(N;k,v)$ and  explore the equivalence
between optimum CDAs and a special class of COAs. We first establish a benchmark to measure the optimality of
$(d,t)$-CDAs.  The following result can be obtained by employing a similar proof as that of Lemma 2.2 in \cite{STY2012}. We simply
replace  $t$-way interactions by consecutive $t$-way interactions and use the fact stated in Lemma \ref{(d,t)-DA-(d-1,t)-DA}. Thus, the proof is omitted here.

%\begin{lemma} \label{lbound-d=1}
%Suppose that $A$ is a $(1,t)$-CDA$(N;k,v)$ with $t< k$. Then,  $|\rho (A,T)|\ge 2$  for any consecutive $t$-way interaction $T$.
%\end{lemma}

%\pf As a CDA is a special type of CCA, we have $|\rho(A,T)| \ge 1$. Thus, it suffices to show $|\rho(A,T)| \not= 1$
%for any consecutive $t$-way interaction $T$. If not,  suppose that $(x_1,x_2, \ldots,x_k)$ is the unique row of $A$
%that covers $T$. This row  also covers $(k-t)$ consecutive $t$-way interactions other than $T$. It follows that there would be
%at least one  $t$-way interaction $T'$ ($\not= T$) such that $\rho(A,T) \subseteq \rho(A,T')$ under the assumption $t<k$.
%Thus, $A$ is not a $(1,t)$-CDA$(N;k,v)$.\qed

%The following lemma can be viewed as a generalization of Lemma \ref {lbound-d=1}.  The proof is similar to that for .

\begin{lemma} \label{lbound-d>1}
Suppose that $A$ is a $(d,t)$-CDA$(N; k,v)$ with $t < k$. Then, $|\rho (A,T)|\geq d+1 $ for any consecutive $t$-way interaction $T$.
\end{lemma}

By applying Lemma \ref {lbound-d>1}, we have a  lower bound on the function $(d,t)$-CDAN$(k,v)$. This serves as our benchmark for measuring
the optimality of CDAs.

\begin{theorem}\label{Lbound-G}
Let $t, k$, and $v$ be positive integers with $t< k$. Then,
\begin{align*}
(d,t)\ \mbox{-CDAN}\ (k,v) \geq & (d+1)v^t.
 \end{align*}
\end{theorem}

\pf Let $A$ be a $(d,t)$-CDA$(N; k,v)$ over $V$ with  $t < k$ and $N = (d,t)$-CDAN$(k,v)$. By definition, for the first $t$ columns $\{1, 2, \cdots, t\}$,
there exist exactly $v^t$ $t$-way interactions of $A$: $\{(i, x_i): 1\leq i\leq t, x_i\in V \}$.
From Lemma \ref {lbound-d>1}, $|\rho (A,T)|\geq d+1 $ for any consecutive $t$-way interaction $T$ of $A$.  Therefore,
$A$ must contain at least $(d + 1)v^t$ rows. This means that $N=(d,t)$-CDAN$(k,v) \geq (d+1)v^t$. \qed

We call a $(d,t)$-CDA$(N;k,v)$ with $N=(d+1)v^t$ {\em optimum}. It is interesting that optimum CDAs have useful applications  in software testing, because they contain minimum rows. In addition, optimum CDAs can be characterized in terms of a special class of COAs. To explore the combinatorial features of optimum CDAs, we need to introduce the notion of simple COAs. A COA$_{\lambda}(t, k, v)$ is  {\em simple} if any $N\times (2t-i)$ subarray consisting of two consecutive $t$ columns with $i$ columns in common contains each $(2t-i)$-tuple at most once, where $0\leq i\leq t-1$. From this definition, it is obvious that a simple COA$_{\lambda}(t, k, v)$ can only exist if $\lambda \le v$. In fact,  a simple COA$_{\lambda}(t, k, v)$ with $\lambda=1$ is a COA$(t,k,v)$.

\begin{example}\label{Ex2-1} The transpose of the following
array is a simple COA$_3(2, 6, 3)$ over $\Z_3.$
 \begin{center}
{\small
 \tabcolsep 2.0pt
 \begin{tabular}{|ccccccccccccccccccccccccccc|}
 \hline
0 & 0 & 0 & 0 & 0 &  0 & 0 & 0 & 0 & 1 & 1 & 1 & 1 & 1 & 1 & 1 & 1 & 1 & 2 & 2 & 2 & 2 & 2 & 2 & 2 & 2 & 2 \\
0 & 0 & 0 & 1 & 1 &  1 & 2 & 2 & 2 & 0 & 0 & 0 & 1 & 1 & 1 & 2 & 2 & 2 & 0 & 0 & 0 & 1 & 1 & 1 & 2 & 2 & 2 \\
0 & 1 & 2 & 0 & 1 &  2 & 0 & 1 & 2 & 0 & 1 & 2 & 0 & 1 & 2 & 0 & 1 & 2 & 0 & 1 & 2 & 0 & 1 & 2 & 0 & 1 & 2 \\
0 & 2 & 1 & 2 & 1 &  0 & 1 & 0 & 2 & 2 & 1 & 0 & 1 & 0 & 2 & 0 & 2 &
1 & 1 & 0 & 2 & 0 & 2 & 1 & 2 & 1 & 0\\
0 & 0 & 0 & 0 & 0 &  0 & 0 & 0 & 0 & 1 & 1 & 1 & 1 & 1 & 1 & 1 & 1 & 1 & 2 & 2 & 2 & 2 & 2 & 2 & 2 & 2 & 2 \\
0 & 1 & 2 & 0 & 1 &  2 & 0 & 1 & 2 & 0 & 1 & 2 & 0 & 1 & 2 & 0 & 1 & 2 & 0 & 1 & 2 & 0 & 1 & 2 & 0 & 1 & 2 \\
\hline
\end{tabular}
}
\end{center}
\end{example}

It is easy to check that any two consecutive columns contain each pair exactly three times, and each 4-tuple over $\Z_3$ from the two disjoint
consecutive two-columns occurs at most once.  For any two consecutive two-columns with one column in common, each 3-tuple occurs at most once. \qed
%\mbox{\quad}

The following theorem explores the equivalence between optimum CDAs and simple COAs.

\begin{theorem} \label{COA-DTA}
Suppose that  $t$ and $k$ are two positive integers and $t< k$.
Then, a simple COA$_{d+1}(t, k, v)$ is equivalent to an optimum $(d,t)$-CDA$((d+1)v^t;k,v)$.
\end{theorem}

\pf ($\Rightarrow$) Suppose that $A$ is a simple COA$_{d+1}(t, k, v)$.
Let  $T$ be an arbitrary consecutive $t$-way interaction of $A$. Consider the set ${\cal T} =\{T_1,T_2,\ldots, T_d\}$
of arbitrary consecutive $t$-way interactions of cardinality $d$  such that $T\not \in {\cal T}$.  We only need to prove that
$\rho(A,T) \not \subseteq \rho(A, \cal {T})$. As $A$ is a COA$_{d+1}(t, k, v)$,  $|\rho(A,T)| = d + 1$ holds for any
consecutive $t$-way interaction. If $\rho(A,T)  \subseteq \rho(A, \cal {T})$, then
there would be at least one $T_j \in \cal {T}$ such that $|\rho(A,T_j)\cap\rho(A,T)|\geq 2$. Suppose that the column indices of $T$ and
$T_j$ have $i$ columns in common, where $0\leq i\leq t-1$. Let $T'$ be the interaction given by deleting the common elements of $T$ and $T_j$ from $T$.
Because $T\not = T_j$, at least two rows in a certain $(2t-i)$ columns would cover the interaction given by a concatenation of $T_j$  and $T'$.
This is inconsistent with the simple property of $A$.

 ($\Leftarrow$) Let $B$ be a $(d,t)$-CDA$((d+1)v^t;k, v)$ over the symbol set $V$.
From Lemma \ref{lbound-d>1}, we know that  $|\rho(A,T)| = d+1$ for any consecutive $t$-way interaction $T$, because
$B$ contains precisely $(d+1)v^t$ rows.  This indicates that each $t$-tuple  occurs as a row exactly $(d+1)$ times in any
consecutive $t$ columns of $B$. Thus, $B$ is a COA$_{d+1}(t,k,v)$. The simple property of $B$ can easily be obtained from the definition of a CDA.\qed

We remark that, by Lemma \ref{NC-CDA}, a $(d,t)$-CDA$(N; k,v)$ cannot exist whenever $d\geq v$.  This is consistent with the fact that a simple COA$_{d+1}(t,k,v)$ can only exist if
$d+1\leq v$.

\section{Construction and Existence of Optimum CDAs }

In this section, we use the equivalence characterization shown
in Theorem \ref{COA-DTA} to construct a great number of optimum CDAs in terms of simple COAs. First, we describe the notion of
super-simple orthogonal arrays (SSOA), which is analogous to the notion of simple COAs.  We say that an OA$_{\lambda}(t, k, v)$ is  {\em
super-simple} if any $(t+1)$ columns of the array contain every $(t+1)$-tuple of symbols as a row at most once. This is denoted by SSOA$_\lambda(t,k,v)$.
Clearly, an  SSOA$_\lambda(t,k,v)$ is a simple COA$_\lambda(t,k,v)$, but the converse is not always true.  It is remarkable that   construction and existence of simple COA$_\lambda(t,k,v)$ with $\lambda>1$ are
considered unless otherwise specified, as we only treat a $(d,t)$-CDA with $d>0$ in this paper.

\subsection{Optimum CDAs from Orthogonal Arrays}

OAs are a highly structured family of arrays that were first introduced by Rao \cite{R1947}. They are also
important objects in combinatorics and experimental design. Over the past half a century, OAs have been the subject of considerable
study. The following elegant result on the existence of OA$(t,k,v)$ with $t\geq 3$ was derived by Bush \cite{Bushb}.

\begin{lemma}\label{Bush1}{\rm\cite{Bushb}}
If $q$ is a prime power and $t< q$, then an OA$(t,q+1,q)$ exists. Moreover, if $q\geq 4$ is a power of 2, an OA$(3,q+2,q)$ exists.
\end{lemma}

Bush also established a powerful composite construction that serves to obtain  new OAs from old ones.
The derived array is formed by juxtaposing these known arrays.

\begin{lemma}\label{Bush2}{\rm \cite{Busha}}
If OA$(t,k,v_i)'$ exists for all $1\leq i\leq m$, then an OA$(t,k,\prod_{i=1}^m v_i)$ exists.
\end{lemma}

Lemmas \ref{Bush1} and \ref{Bush2} immediately give the following result.

\begin{lemma}\label{Bush}
Suppose that $v=q_1q_2\cdots q_s$ is a standard factorization of $v$ into distinct prime powers.
If $q_i>t$, then an OA$(t,k+1,v)$ exists, where $k=\mbox{min}\{q_i:1\leq i\leq s\}$.
\end{lemma}

The following well-known result (see, for
example, \cite{CD2007,HSS1999}) can be obtained by zero-sum construction, i.e., for each of the $v^t$ $t$-tuples over $\Z_v$,
form a row vector of length $t+1$ by adjoining the negative of the sum of the elements in the first $t$ columns to the last column.

\begin{lemma}\label{OA(t,t+1,v)}
An OA$(t,t+1,v)$ exists for any integer $v\geq 2$ and $t\geq 2$.
\end{lemma}

The existence  of OAs with $t=3$ and $k=5,6$ was recently proved by Ji
and Yin \cite{JY2010}.

\begin{lemma}\label{OA(3,5,v)}{\rm\cite{JY2010}}
Let $v\geq 4$ be an integer. If $v\not\equiv2 \pmod 4$, then an
OA$(3,5,v)$ exists.
\end{lemma}

\begin{lemma}\label{OA(3,6,v)}{\rm\cite{JY2010}}
Let $v$ be a positive integer that satisfies gcd$(v,4)\not=2$ and
gcd$(v,18)\not=3$. Then, there is an OA$(3,6,v)$, and
OA$(3,6,3u)$ with $u\in \{5,7\}$ exists.
\end{lemma}

Our approach of constructing optimum CDAs from SSOAs is as follows.

\begin{construction}\label{SSOA-COA}
If an SSOA$_\lambda(t,k,v)$ exists, then a simple COA$_{\lambda}(t,k+t-1,v)$ exists.
\end{construction}
\pf Let $A$ be the given SSOA$_\lambda(t,k,v)$ over the symbol set $V$ with column vectors $A_1,A_2,\cdots,A_k$. Write
$A'=(A_1,A_2,\cdots, A_k,A_1,A_2,\cdots, A_{t-1})$. We claim that the array $A'$ is a simple COA$_{\lambda}(t,k+t-1,v)$ over $V$, as
desired.  Clearly, it is a COA$_\lambda(t,k,v)$ because any consecutive $t$ columns  are the certain $t$ columns of $A$. It remains to prove that $A'$
is simple. For any two consecutive $t$ columns $i,i+1,\cdots,i+t-1$ and $j,j+1,\cdots,j+t-1$, suppose that $|\{i,i+1,\cdots,i+t-1\}\cap \{j,j+1,\cdots,j+t-1\}|=l$ with
$0\leq l\leq t-1$, where $1\leq i\leq k-1, i<j\leq k$. If $l=t-1$, then the consecutive $t+1$ columns are the certain $t+1$ columns of $A$. As $A$ is an SSOA$_\lambda(t,k,v)$,
each $(t+1)$-tuple occurs at most once. Hence, any $N\times (t+1)$ subarray consisting of the consecutive $t+1$ columns in $A'$ contains each $(t+1)$-tuple at most once.
This can be divided into two cases when $l=0$.  If two disjoint consecutive $t$ columns are from the first $k$ columns of $A'$, each $2t$-tuple occurs at most once because each $(t+1)$-tuple
occurs at most once in $A$. If two disjoint consecutive $t$ columns do not lie in the first $k$ columns of $A'$ at the same time, then there are at least $(t+1)$ columns in the $2t$ columns that lie in the columns
of $A$. The super-simple property of $A$ guarantees that each $(t+1)$-tuple occurs at most once, and thus any $N\times 2t$ subarray contains each $2t$-tuple at most once. Similarly,
if  $1\leq l\leq t-2$,  we can prove that any $N\times (2t-l)$ subarray is simple.
\mbox{\quad} \qed

As mentioned in  Construction 3.7  \cite {STY2012}, let $A$ be  an  OA$(t + 1, k + 1, v)$ over symbol set $V$.
Taking the rows in $A$ that begin with $\lambda$ distinct symbols $x_i\ (i=1, 2, \cdots, \lambda)$ from $V$ and
omitting the first column yields an SSOA$_\lambda(t, k, v)$. By Construction \ref{SSOA-COA}, we have the following working
method of construction.

\begin{construction}\label{OA-COA}
If an OA$(t+1,k+1,v)$ exists, then a simple COA$_{\lambda}(t,k+t-1,v)$ exists for any positive integer $ \lambda \leq v$.
\end{construction}

The following example illustrates the idea of Construction \ref {OA-COA}.

\begin{example}\label {Ex3-1}\rm
Let $A=\left(
\begin{array}{c}
A_0^T\\
A_1^T\\
A_2^T
\end{array}
\right)$,  where
{\small
$
\begin{array}{lll}
A_0=
 \tabcolsep 0.5pt
 \begin{tabular}{|ccccccccc|}
 \hline
0 &0 &0 &0 &0 &0 &0 &0 &0 \\
0 &0 &0 &1 &1 &1 &2 &2 &2 \\
0 &1 &2 &0 &1 &2 &0 &1 &2 \\
0 &2 &1 &2 &1 &0 &1 &0 &2 \\
\hline
\end{tabular}
&
A_1=
\tabcolsep 0.8pt
\begin{tabular}{|ccccccccc|}
\hline
1 &1 &1 &1 &1 &1 &1 &1 &1 \\
0 &0 &0 &1 &1 &1 &2 &2 &2 \\
0 &1 &2 &0 &1 &2 &0 &1 &2 \\
2 &1 &0 &1 &0 &2 &0 &2 &1 \\
\hline
\end{tabular}
&
A_2=
\tabcolsep 0.8pt
\begin{tabular}{|ccccccccc|}
\hline
2 &2 &2 &2 &2 &2 &2 &2 &2 \\
0 &0 &0 &1 &1 &1 &2 &2 &2 \\
0 &1 &2 &0 &1 &2 &0 &1 &2 \\
1 &0 &2 &0 &2 &1 &2 &1 &0 \\
\hline
\end{tabular}
\end{array}
$}

\noindent
Then $A$ is an OA$(3,4,3)$ over $\Z_3$. We  take the rows which involve the symbols $0$ and $2$ in the first column of $A$ to form an SSOA$_2(2,3,3)$ $A'$, where

{\small
$
\begin{array}{lll}
A'=\left(
\begin{array}{c}
A_0'^T\\
A_2'^T
\end{array}
\right),&
A_0'=
 \tabcolsep 0.5pt
 \begin{tabular}{|ccccccccc|}
 \hline
0 &0 &0 &1 &1 &1 &2 &2 &2 \\
0 &1 &2 &0 &1 &2 &0 &1 &2 \\
0 &2 &1 &2 &1 &0 &1 &0 &2 \\
\hline
\end{tabular}
&
A_2'=
\tabcolsep 0.8pt
\begin{tabular}{|ccccccccc|}
\hline
0 &0 &0 &1 &1 &1 &2 &2 &2 \\
0 &1 &2 &0 &1 &2 &0 &1 &2 \\
1 &0 &2 &0 &2 &1 &2 &1 &0 \\
\hline
\end{tabular}
\end{array}
$}

\noindent

Let $T=(00 0 1 1 1 2 2 200 0 1 1 1 2 2 2)^T$. By Construction \ref{SSOA-COA}, the array $(A'|T)$ is a simple COA$_2(2,4,3)$ over $\Z_3$. \qed

\end{example}

If an OA$(3,k+1,v)$ exists, then a simple COA$_\lambda(2,k+1,v)$ exists by Construction \ref{OA-COA}. In fact, we can improve
the number of factors as follows.

\begin{construction}\label{OA3-COA2}
If an OA$(3,k+1,v)$ with $k>4$ exists, then a simple COA$_\lambda(2,2k+1,v)$ exists for any integer $\lambda \leq v$.
\end{construction}

\pf If an OA$(3,k+1,v)$ exists, then an SSOA$_\lambda(2,k,v)$ $(A_i),   i\in I_k$ exists for any integer $\lambda \leq v$.

Write
\begin{eqnarray*}
A' = \left\{
\begin{array}{ll}
(A_1,A_2,\cdots, A_k,A_1,A_3,\cdots, A_{k-2}, A_k, A_2,A_4,\cdots,A_{k-3}, A_{k-1}, A_1) & \mbox{if}\  k \ \mbox{is odd} ;\\
(A_1,A_2,\cdots, A_k,A_1,A_3,\cdots, A_{k-3}, A_{k-1}, A_2, A_4,\cdots,A_{k-2}, A_k, A_2), & \mbox{if}\  k \ \mbox{is even}.
\end{array}
\right.
\end{eqnarray*}

\noindent
It is easily checked that $A'$ is the required simple COA. \qed

Combining Theorem \ref{COA-DTA} and Constructions \ref{OA-COA} and \ref{OA3-COA2} with  those known OAs given in the previous
lemmas, it is possible to produce an infinite series of optimum CDAs.

\begin{theorem}\label{COA(t,2t,v)}
Let both $t\ge 2$ and $v\ge 2$ be integers. Then, an optimum  $(d,t)$-CDA$((d+1)v^t;2t,v)$ exists for any positive integer $d$ satisfying $d+1\le v$.
\end{theorem}

\pf Under the given assumption, an OA$(t+1, t+2, v)$ exists by Lemma \ref {OA(t,t+1,v)}.  Applying Construction \ref{OA-COA} with this OA will produce a simple COA$_{d+1}(t, 2t, v)$, because $d+1\le v$.  The conclusion follows from Theorem \ref{COA-DTA} with $k=2t$. \qed

\begin{theorem}
Let $q$ be a prime power. Then, an optimum $(d,t)$-CDA$((d+1)q^t;t,q+t-1,q)$ exists for any positive integers $d$ and $t+1$
that are less than $q$.  Moreover, if $q\geq 4$ is a power of 2,  then an optimum $(d,2)$-CDA$((d+1)q^2;2q+3,q)$ also exists.
\end{theorem}

\pf  Apply Theorem \ref{COA-DTA}. The required simple COAs are
obtained  by Constructions \ref{OA-COA} and \ref{OA3-COA2} and Lemma \ref{Bush1}. \qed

\begin{theorem}
Suppose that $v=q_1q_2\cdots q_s$ is a standard factorization of $v$ into distinct prime powers and $k=\mbox{min}\{q_i:1\leq i\leq s\}$.
If $q_i>t+1$ and $d+1\le v$, then there exists an optimum $(d,t)$-CDA$((d+1)v^t;k+t-1,v)$.
\end{theorem}

\pf By Lemma \ref{Bush}, we have an OA$(t+1, k+1, v)$. Applying
Construction \ref{OA-COA} produces a simple COA$_{d+1}(t, k+t-1, v)$ under the assumption $d+1\le v$.
Hence, we can  apply  Theorem \ref{COA-DTA}  to  form  the desired  CDA. \qed

\begin{theorem}
Let $v$ be a positive integer satisfying gcd$(v,4)\not=2$ and gcd$(v,18)\not=3$. Then, there
exists an optimum $(d,2)$-CDA$((d+1)v^2;11,v)$ for any  positive integer $d$ with $d+1\le v$.
\end{theorem}

\pf Apply Theorem \ref{COA-DTA}. The required simple  COAs are
given by Construction \ref{OA3-COA2} and Lemma \ref{OA(3,6,v)}. \mbox{\quad} \qed

\begin{theorem}\label{COA(2,6,v)}
Let $v\geq 4$ be an integer. If $v\not\equiv2 \pmod 4$, then an optimum $(d,2)$-CDA$((d+1)v^2;6,v)$ exists for any positive integer $d$ with $d+1\le v$.
\end{theorem}

\pf Under the given assumption, an OA$(3,5,v)$ exists by Lemma \ref{OA(3,5,v)}. It can be used to produce an SSOA$((d+1)v^2;4,v)$, $A$ by Construction 3.7 in
\cite{STY2012}, where $d+1\le v$. Let $A_1,A_2,A_3,A_4$ be the column vectors of $A$. The array $(A_1,A_2,A_3,A_4,A_1,A_3)$ is the desired optimum CDA. \qed

\subsection{Some more approaches for constructing  simple COAs }

In this subsection, we present some more constructions of optimum CDAs  in design theory.
Let $A$  be a COA$_{\lambda}(t, k, v)$ over the symbol set $V$. If the rows of $A$ can be partitioned into
$\mu$ subarrays  such that  each has the simple property described above, then  $A$ is termed a {\em
$\mu$-row-divisible COA$_{\lambda}(t, k, v)$}. Two simple COAs over the same symbol set are {\em compatible}
if their superimposition constitutes a simple COA. A set of $w$ simple COAs over the same symbol set is termed
{\em compatible} if all elements of the set are pairwise compatible. The notion of $\mu$-row-divisible and compatible OAs was first introduced
in \cite{STY2012}, where they were used to construct SSOAs. Here, we modify them to construct simple COAs.

\begin{theorem} \label{GCon}
Let $v_1$ and $v_2$ be two positive integers such that $\mu$ compatible simple COA$_{\eta}(t, k, v_2)'$s exist.
Suppose that there are $r$ non-negative integers $m_1,m_2,\cdots,m_r$ and $2r$ positive integers $\mu_1, \mu_2, \cdots, \mu_r$, $
\lambda_1,\lambda_2,\cdots, \lambda_r$ such that the following two conditions are both satisfied:
\begin{enumerate}
\item $m_1\mu_1+m_2\mu_2+\cdots+m_r\mu_r\ \leq \mu$;
\item a $\mu_i$-row-divisible COA$_{\lambda_i}(t, k, v_1)$ exists for $1\leq i\leq r$.
\end {enumerate}
\noindent Then, there exists a simple COA$_{\eta(m_1\lambda_1+m_2\lambda_2+\cdots+m_r\lambda_r)}(t, k,
v_1v_2)$.
\end{theorem}

\pf Without loss of generality, we can assume that all the given row-divisible COAs are defined
on the same symbol set $V$ (otherwise, we may take an appropriate permutation of the symbols). For each $i$ with $1\leq i\le r$,
we form an $m_i\mu_i$-row-divisible COA$_{m_i\lambda_i}(t, k, v_1)$  by taking $m_i$ copies of  a $\mu_i$-row-divisible COA$_{\lambda_i}(t, k, v_1)$.
The superimposition of the resultant $r$ row-divisible OAs is then the $(m_1\mu_1+m_2\mu_2+\cdots
+m_r\mu_r)$-row-divisible OA$_{m_1\lambda_1+m_2\lambda_2+\cdots+m_r\lambda_r}(t,k,v_1)$, as desired.

By assumption, there exist $\mu $ compatible simple COA$_{\eta}(t, k, v_2)'$s. Thus, the desired simple COAs
can be obtained by  a modification of the usual weighting method in
design theory, as used in the proof of \cite [Lemma 13]{H2000} and
\cite [Construction 4.2]{STY2012}. \mbox{\quad}\qed

Here, we give an example for the superimposition of row-divisible COAs. A COA$_2(3,5,2)$ is given by Lemma \ref{simple COA(2,5,236)} and denoted as $A$.
A $2$-row-divisible COA$_3(2,5,2)$ with two partitions $B_1$ and $B_2$  is given by Lemma \ref{2-row-divisible COA(2,5,2)}. Then $\left(
\begin{array}{c}
A\\
B_1\\
B_2
\end{array}
\right)$ is a $3$-row-divisible COA$_5(2,5,2)$ with $3$ partitions $A$, $B_1$ and $B_2$. Similarly, $\left(
\begin{array}{c}
A\\
A\\
B_1\\
B_2
\end{array}
\right)$ is a $4$-row-divisible COA$_7(2,5,2)$ with $4$ partitions $A$, $A$, $B_1$ and $B_2$.

By taking $r=1$, $m_1=1$, $\lambda_1=\lambda$, $v_1=v$, $\mu_1=\mu$, and $v_2=m$ in Theorem \ref {GCon}, we
obtain the following corollary.

\begin {corollary} \label{Inflation-G}
Let $v$,  $k$, and $t$ be integers satisfying $k\ge t\ge 2.$ If a
$\mu$-row-divisible COA$_{\lambda}(t, k, v)$  and $\mu$ compatible
simple  COA$_{\eta}(t, k, m)'$s all exist, then so does a simple COA$_{\lambda \eta}(t, k, mv)$. In particular,  if a simple
COA$_{\lambda}(t, k, v)$ and a simple COA$_{\eta}(t, k, m)$ both
exist, then so does a simple COA$_{\lambda \eta}(t, k, mv).$
\end {corollary}

%By taking $\lambda_1=\lambda$ and $\lambda_2=\eta$ in Corollary \ref{Inflation-G}, we obtain
%the following  corollary.
%\begin {corollary} \label{Inflation-1} Let $v$,  $k$,
%and $t$ be integers satisfying $k\ge t\ge 2$. If a
%$\mu$-row-divisible COA$_{\lambda}(t, k, v)$  and $\mu$ compatible
%COA$(t, k, m)$'s all exist, then so does a simple COA$_{\lambda}(t, k, mv)$.
%\end {corollary}
%
%
%Further, by taking $\mu=1$ in Corollary
%\ref{Inflation-1}, we have the following.
%
%\begin {corollary}  \label{Inflation-2}
%Suppose that a simple COA$_{\lambda}(t,k,v)$ and an OA$(t, k,
%m)$ exist. Then, a simple COA$_{\lambda}(t, k, mv)$ exists.
%\end {corollary}

\subsection{Existence spectrum of optimum CDAs with few factors for $t\in \{2,3\}$}

In this subsection, the existence of optimum CDAs with few factors is determined completely when $t\in \{2,3\}$. It is known that
the derived array of an OA$(t,k,v)$ is an OA$(t-1,k-1,v)$. This simple fact is not always true for consecutive orthogonal arrays. The following example
indicates this fact.

\begin{example}\label{CCOA$(4,6,2)}
The transpose of the following array is a COA$(4,6,2)$ over $\Z_2$.
 \begin{center}
{\small
 \tabcolsep 2.0pt
 \begin{tabular}{|cccccccc|cccccccc|}
 \hline
 0 &  0&   0&   0&   0&   0&   0&   0&   1&   1&   1&   1&   1&   1&   1&   1\\
 \hline
 0 &  0&   0&   0&   1&   1&   1&   1&   0&   0&   0&   0&   1&   1&   1&   1\\
 0 &  0&   1&   1&   0&   0&   1&   1&   0&   0&   1&   1&   0&   0&   1&   1\\
 0 &  1&   0&   1&   0&   1&   0&   1&   0&   1&   0&   1&   0&   1&   0&   1\\
 0 &  0&   1&   1&   0&   0&   1&   1&   1&   1&   0&   0&   1&   1&   0&   0\\
 0 &  1&   0&   1&   1&   0&   1&   0&   0&   1&   0&   1&   1&   0&   1&   0 \\

\hline
\end{tabular}
}
\end{center}

\noindent The arrays $A_0$ and $A_1$ derived by deleting the last column are as follows.

{\small
$$
\begin{array}{ll}
A_0=
 \tabcolsep 1.8pt
 \begin{tabular}{|cccccccc|}
 \hline
 0 &  0&   0&   0&   1&   1&   1&   1 \\
 0 &  0&   1&   1&   0&   0&   1&   1 \\
 0 &  1&   0&   1&   0&   1&   0&   1 \\
 0 &  0&   1&   1&   0&   0&   1&   1 \\
 0 &  1&   0&   1&   1&   0&   1&   0 \\
\hline
\end{tabular}
\quad
&
\quad
A_1=
\tabcolsep 1.8pt
 \begin{tabular}{|cccccccc|}
 \hline
0&   0&   0&   0&   1&   1&   1&   1 \\
0&   0&   1&   1&   0&   0&   1&   1 \\
0&   1&   0&   1&   0&   1&   0&   1 \\
1&   1&   0&   0&   1&   1&   0&   0 \\
0&   1&   0&   1&   1&   0&   1&   0 \\
\hline
\end{tabular}
\end{array}
$$
}
\noindent  Clearly, neither $A_0$ nor $A_1$ is a COA$(3,5,2)$. \qed
\end{example}

Example \ref{CCOA$(4,6,2)} tells us that the existence of a COA$(t,k,v)$ does not imply the existence of a COA$(t-1,k-1,v)$.
Moreover, deleting a certain column from a COA$(t,k,v)$  does not guarantee it is a COA$(t,k-1,v)$.  Thus, we have to construct simple COAs for each set of values $t,k,v$.

\begin{theorem}
An optimum $(d,2)$-CDA$(k,v)$ with $k=3,4$ exists for any integer $d$ with $(d+1)\leq v$.
\end{theorem}

\pf The existence of SSOA$_{d+1}(2,3,v)$ is proved in \cite{STY2012}. Clearly, an SSOA is also a
simple COA. A simple COA$_{d+1}(2,4,v)$ can be obtained by taking $t=2$ in Theorem \ref{COA(t,2t,v)}. Applying Theorem \ref{COA-DTA} produces
optimum CDAs, as desired. \qed

Similarly, we can obtain the following result.

\begin{theorem}
An optimum $(d,3)$-CDA$(4,v)$  exists for any integer $d$ with $(d+1)\leq v$.
\end{theorem}

The following theorem treats the case with $t=3,k=5,6$.
\begin{theorem}
An optimum $(d,3)$-CDA$(k,v)$  with $k=5,6$ exists for any integer $d$ with $(d+1)\leq v$.
\end{theorem}

\pf From Theorem \ref{COA-DTA}, we only need to construct a simple COA$_{d+1}(3,5,v)$ and COA$_{d+1}(3,6,v)$ for any integer $d$ with $(d+1)\leq v$.
We take $t=3$ in Construction \ref{OA-COA} with an OA$(4,5,v)$ to form  arrays $A'=(A_1,A_2,A_3,A_4,A_1)$ and $A''=(A_1,A_2,A_3,A_4,A_1,A_2)$,  where $A_i$ is
the $i$th column of $A$. It is routine to check that $A'$ and $A''$ are the simple COA$_{d+1}(3,5,v)$ and COA$_{d+1}(3,6,v)$, respectively. \mbox{\quad} \qed

By taking the array $(A_1,A_2,A_3,A_4,A_1)$ in Theorem \ref{COA(2,6,v)}, we obtain the following result.

\begin{theorem}\label{vno=4t+2}
Let $v\geq 4$  be an integer.  If $v\not\equiv2 \pmod 4$, then an optimum $(d,2)$-CDA$(5,v)$  exists for any integer $d$ with $(d+1)\leq v$.
\end{theorem}

For the completeness  of existence for a $(d,2)$-CDA$(5,v)$, we consider the case $v=2,3$ or $v\geq 6$ and $v \equiv2 \pmod 4$. For $v\in \{2,3,6\}$, we have the following results.

\begin{lemma}\label{2-row-divisible COA(2,5,2)}
A $2$-row-divisible  COA$_3(2,5,2)$ over $\Z_2$ exists.
\end{lemma}
\pf Simply take the transpose of the following array with two partitions:

\begin{center}
\tabcolsep 2.8pt
 \begin{tabular}{|cccccc|cccccc|}
 \hline
0&   0&   1&   1&   0&   1&   0&   1&   0&   1&   1&   0    \\
0&   0&   1&   0&   1&   1&   1&   0&   1&   1&   0&   0    \\
0&   1&   0&   1&   1&   1&   0&   1&   1&   0&   0&   0    \\
0&   0&   1&   1&   1&   0&   1&   1&   0&   0&   0&   1    \\
0&   1&   1&   1&   0&   0&   1&   0&   1&   0&   1&   0    \\
\hline
\end{tabular}
\end{center}
\qed

\begin{lemma}\label{simple COA(2,5,236)}
A simple COA$_v(2,5,v)$ over $\Z_v$ exists for any integer $v\geq 2$.
\end{lemma}
\pf By Lemma \ref{OA(t,t+1,v)}, an OA$(3,4,v)$  with column vectors  $A_1,A_2, A_3,A_4$ exists.
Write $A'=(A_1,A_2,\\A_3,A_4,A_1)$.  $A$ is the required simple COA. \qed

\begin{lemma}\label{simple COA(2,5,23)}
A simple COA$_\lambda(2,5,v)$ over $\Z_v$ exists for $(\lambda,v)\in \{(2,3),(1,6),(2,6),(3,6),(4,6),(5,6)\}$.
\end{lemma}

\pf The existence of an SSOA$_2(2,5,v)$ for $v=3,6$ is proved in \cite{Chen2011}. This implies the existence of simple
COA$_2(2,5,3)$ and COA$_2(2,5,6)$.  A simple COA$_4(2,5,6)$ can be obtained  by Corollary \ref{Inflation-G}. The ingredient
COA$_2(2,5,2)$ is given by Lemma \ref{simple COA(2,5,236)}. A COA$(2,5,6)$ is constructed using the array $(A_1,A_2,A_3,A_1,A_3)$,  where $A_i$ is
the $i$th column of an OA$(2,3,6)$. A COA$_\lambda(2,5,6)$ with $\lambda=3$ or 5 is given in the Appendix.\qed

%\begin{lemma}\label{2-row-divisible COA9(2,5,6)}
%A $2$-row-divisible  COA$_\lambda(2,5,6)$ over $\Z_6$ with $\lambda\in \{3,5,9\}$ exists.
%\end{lemma}
%
%\pf Let $A$ and $B$ be a   $2$-row-divisible COA$_3(2, 5, 2)$ with the partition $A_1, A_2$  and a simple COA$_3(2,\\5,3)$, respectively. Take $V=\Z_2\times Z_3$. Over
%$V$, we form a $9\cdot 6^2\times 5$ array $\overline{A}$ as follows. For each row $(a_{i1}, a_{i2 }, \cdots,
%a_{i5})$ of $A_i$ and each row $(b_{h1}, b_{h2},\cdots , b_{hk})$ of
%$B$ for $i=1,2$, include the row $\{(a_{i1}, b_{h1}), (a_{i2}, b_{h2}), \cdots, (a_{ik},
%b_{hk})\}$ in $\overline{A}$ as a row. It is easy to check that the resultant array $\overline{A}$ is
%a $2$-row-divisible COA$_9(2,5,6)$. A $2$-row-divisible  COA$_\lambda(2,5,6)$ over $\Z_6$ with $\lambda=3,5$ can be obtained by
%the juxtaposition of an OA$(2,5,6)$ and a simple COA$_2(2,5,6)$ or a simple COA$_4(2,5,6)$, respectively. \mbox{\quad}\qed

We now determine  the existence for a simple COA$_\lambda(2,5,v)$ with $v\equiv 2 \pmod 4$. Write $v=4t+2=2(2t+1)$, where $t\geq 2$.
To apply the corollary in Subsection 4.2, we need compatible COAs, which can be obtained using the simple argument in \cite{STY2012}.

\begin{lemma}\label{compatible COA(2,5,v)}
Let $v=2t+1$ be an integer with  $t\geq 2$. If $2t+1\not =3u$ with gcd$(u,6)=1$, then $v$ compatible COA$(2,5,v)'$s exist.
Moreover, if $v=2t+1=3u$ with $u\not =1$ and gcd$(u,6)=1$, the $u$ compatible COA$(2,5,u)'$s exist.
\end{lemma}
\pf Under the given assumption and Lemma \ref{OA(3,6,v)}, an OA$(3,6,v)$ exists. By Construction 3.7 in \cite{STY2012}, the
$v$ derived arrays given by deleting the first column form $v$ compatible COAs. The second assertion can be proved in a similar way to the first statement. \qed

\begin{lemma}\label{simple COA(2,5,4t+2)}
Let $v=4t+2$ be an integer with $t\geq 2$. If $2t+1\not =3u$ with gcd$(u,6)=1$, then a simple COA$_\lambda(2,5,v)$ exists for any integer $\lambda \leq v$.
\end{lemma}
\proof   From Lemmas \ref{2-row-divisible COA(2,5,2)} and \ref{simple COA(2,5,236)},
a $\mu_i$-row-divisible OA$_{\lambda_i}(2, 5, 2)$  exists for $i=1, 2.$   Apply Theorem \ref {GCon} with $v_1=2$, $v_2=2t+1$, $(\mu_1, \lambda_1)=(1, 2)$,
and $(\mu_2, \lambda_2)=(2, 3)$. It remains to show that the system of equations
$$ \left\{
\begin{array}{l}
2m_1+3m_2=\lambda ,\\
m_1+2m_2\le 2t+1.
\end{array}
\right.
$$
is solvable for non-negative integers $m_1$ and $m_2$
for any given $\lambda$ with $ \lambda \le 4t+2$.

\noindent
It now turns out that
$$
(m_1, m_2) = \left\{
\begin{array}{ll}
(\frac{\lambda}{2}, 0),\ \mbox{if\ $\lambda$ is even },\\

(\frac{\lambda-3}{2}, 1),\ \mbox{if\ $\lambda$ is odd},
\end{array}
\right.
$$
\noindent is one solution of the above system of equations.
\mbox{\quad}\qed

\begin{lemma}\label{simple COA(2,5,3u)}
Let $v=6u$ be an integer with $u\not =1$ and gcd$(u,6)=1$. Then, a simple COA$_\lambda(2,5,v)$ exists for any integer $\lambda \leq v$.
\end{lemma}

\pf From Lemmas \ref{simple COA(2,5,236)} and \ref{simple COA(2,5,23)}, we know that a simple COA$_{\lambda}(2,5, 6)$ with $\lambda \in \{1, 2, 3,4,5, 6\}$ exists. For any given $\lambda$ with $ \lambda \le 6u$, we write $h = \lfloor
\lambda / 6\rfloor$ and $\varepsilon = \lambda - 6h \in \{0,1,\cdots,5\}$. Then,
$$ \left\{
\begin{array}{ll}
h \le u,\ \mbox{if\ $\varepsilon\ =\ 0$},\\
h \le u - 1,\ \mbox{if\ $\varepsilon\ =1,2,3,4,5$}.
\end{array}
\right.
$$

\noindent Apply Theorem \ref {GCon}  with $v_1=6$, $v_2=u$, $(\mu_1, \lambda_1)=(1, 6)$, and $(\mu_2, \lambda_2)=(1,\varepsilon)$, where $\varepsilon=1,2,3,4,5$

It can easily be checked that the system of equations
$$ \left\{
\begin{array}{l}
6m_1+\lambda_2 m_2=\lambda,\\
 m_1+\mu_2 m_2\le u.
\end{array}
\right.
$$
has one solution with non-negative integers:
$$
(m_1, m_2) = \left\{
\begin{array}{ll}
(h, 0),\ \mbox{if\ $\varepsilon\ =\ 0$};\\
(h, 1),\ \mbox{if\ $\varepsilon\ =\ 1, 2$, 3, 4, 5}.\\
\end{array}
\right.
$$
\qed

Combining Theorem \ref{COA-DTA} with the results in Theorem \ref{vno=4t+2}, Lemmas \ref{simple COA(2,5,236)}, \ref{simple COA(2,5,23)}, \ref{simple COA(2,5,4t+2)}, \ref{simple COA(2,5,3u)}, we have the following results.

\begin{theorem}
An optimum $(d,2)$-CDA$(5,v)$ exists for any positive integer $d+1\leq v$.
%xcept possibly where
%
%\begin{enumerate}
%\item $(d,v)\in \{(2,6),(4,6)\}$;
%\item $v=6u$ and $d=v-2$, where $u\not =1$ and gcd$(u,6)=1$.
%\end{enumerate}

\end{theorem}

\section{Concluding Remarks}
Detecting arrays are used to generate test suites for locating and detecting interaction faults between factors.
For practical software testing, there may only be interactions between neighboring factors. Although DAs can be used
to locate and detect interaction faults between neighboring factors,  they are not well adapted for this kind of
software testing.  This paper has introduced the notion of consecutive detecting arrays to solve this problem. Consecutive
detecting arrays are of interest in generating software test suites to cover any consecutive  $t$-way component interactions and locate interaction faults between neighboring factors. In this paper, a general lower bound on the size of $(d,t)$-CDA$(N;k,v)$ has been established. The equivalence between the optimum  $(d,t)$-CDA$((d+1)v^t;k,v)$ and a simple COA$_{d+1}(t, k, v)$
was explored in Theorem \ref{COA-DTA}. Based on this equivalence, a great number of optimum CDAs that satisfy the lower bound were obtained by constructing simple COAs. The existence spectrum  of $(d,t)$-CDA$((d+1)v^t;k,v)$ with few factors for $t=2,3$ was completely determined. Future studies will focus on finding new techniques for constructing simple  COAs and deriving more results with large numbers of factors. In particular, just as consecutive
orthogonal arrays play a role in this respect, we can further study circular consecutive detecting arrays by means of circular consecutive covering arrays proposed in \cite{RMS2018}.

\vspace{0.2cm}
Appendix

\vspace{0.2cm}
A simple COA$_\lambda(2,5,6)$ with $\lambda=3 $ or $5$ is given below.

$\lambda=3$, $A=(A_1,A_2)^T$, where $A_1$ and $A_2$ are as follows.

\vspace{0.2cm}
$A_1=\left(
\begin{array}{c}
000000000000000000111111111111111111222222222222222222\\
000111222333444555000111222333444555000111222333444555\\
013023245145035124245145013023124035013023245145035124\\
453154032210025134201023541345314205501205413324430215\\
111111111111111111222222222222222222333333333333333333
\end{array}
\right)$

\vspace{0.2cm}
$A_2=\left(
\begin{array}{c}
333333333333333333444444444444444444555555555555555555\\
000111222333444555000111222333444555000111222333444555\\
245145013023124035013023245145035124245145013023124035\\
342431052150125340042340521135514320153512403204230154\\
444444444444444444555555555555555555000000000000000000
\end{array}
\right)$

\vspace{0.2cm}
$\lambda=3$, $B=(B_1,B_2,B_3)^T$, where $B_1$, $B_2$ and $B_3$ are as follows.

\vspace{0.2cm}
$B_1=\left(
\begin{array}{c}
000000000000000000000000000000111111111111111111111111111111\\
000001111122222333334444455555000001111122222333334444455555\\
012340123501245013450234512345012450123401345012351234502345\\
012341234523450345014501250123103522145032415435204012351034\\
111111111111111111111111111111222222222222222222222222222222
\end{array}
\right)$

\vspace{0.2cm}
$B_2=\left(
\begin{array}{c}
222222222222222222222222222222333333333333333333333333333333\\
000001111122222333334444455555000001111122222333334444455555\\
013450234501234123450124501235012350124512345023450123401345\\
542313012401354540122315040235421050524305132103245134034251\\
333333333333333333333333333333444444444444444444444444444444
\end{array}
\right)$

\vspace{0.2cm}
$B_3=\left(
\begin{array}{c}
444444444444444444444444444444555555555555555555555555555555\\
000001111122222333334444455555000001111122222333334444455555\\
123450134502345012450123501234023451234501235012340134501245\\
345105023141023213450523412540254034135015024502133254103412\\
555555555555555555555555555555000000000000000000000000000000
\end{array}
\right)$


\begin{thebibliography}{99}


\bibitem{Busha}
K. A. Bush, A generalization of the theorem due to MacNeish, {\it Ann. Math. Stat.}, 23: 293-295, 1952.

\bibitem{Bushb}
K. A. Bush, Orthogonal arrays of index unity, {\it Ann. Math. Stat.}, 23: 426-434, 1952.

\bibitem{CCK1999}
M. Chateauneuf, C. J. Colbourn and D. L. Kreher, Covering arrays of strength three, {\it Des. Codes Cryptogr.}, 16: 235-242, 1999.

\bibitem{CK2002}
M. Chateauneuf and D. L. Kreher, On the state of strength-three covering arrays, {\it J. Combin. Des.},  10: 217-238, 2002.

\bibitem{Chen2011}
Y. Chen, Constructions of Optimal Detecting Arrays of Degree 5 and Strength 2, Master's Thesis, Soochow University, 2011.

\bibitem{Colbourn2004}
C. J. Colbourn, Combinatorial aspects of covering arrays, {\it LeMatematiche (Catania)}, 58: 121-167, 2004.

\bibitem{Colbourn2008}
C. J. Colbourn, Strength two covering arrays: Existence tables and projection, {\it Discrete Math.}, 308: 772-786, 2008.


\bibitem{CD2007}
C. J. Colbourn and J. H. Dinitz, The CRC Handbook of Combinatorial Designs, CRC Press, Boca Raton, FL, 2007.


\bibitem{CMTW2006}
C. J. Colbourn, S. S. Martirosyan, T. V. Trung and R. A. Walker II, Roux-type constructions for covering arrays of strengths three and four, {\it Des. Codes Cryptogr.}, 41: 33-57, 2006.

\bibitem{CM2008}
C. J. Colbourn and D. W. McClary, Locating and detecting arrays for interaction faults, {\it J. Comb. Optim.}, 15: 17-48, 2008.

\bibitem{GKM2010}
A. P. Godbole, M. V. Koutras and F. S. Milienos, Consecutive covering arrays and a new randomness test, {\it J. Statist. Plann. Inference}, 140(5): 1292-1305, 2010.



\bibitem{GKM2011}
A. P. Godbole, M. V. Koutras and F. S. Milienos, Binary consecutive covering arrays, {\it Ann. Inst. Stat. Math.}, 63(3): 559-584, 2011.


\bibitem{H2000}
S. Hartman, On simple and supersimple transversal designs, {\it J. Comb. Des.}, 8: 311-322, 2000.


\bibitem{HR2004}
A. Hartman and L. Raskin, Problems and algorithms for covering arrays, {\it Discrete Math.}, 284: 149-156, 2004.



\bibitem{HSS1999}
A. S. Hedayat,  N. J. A. Slone and J. Stufken, Orthogonal Arrays, Springer, New York, 1999.

\bibitem{JY2010}
L. Ji and J. Yin,  Constructions of new orthogonal arrays and covering arrays of strength three, {\it J. Combin. Theory (A)}, 117: 236-247, 2010.


\bibitem{JM2018}
J. T. Jimenez and  I. I. Marquez,  Covering arrays of strength three from extended permutation vectors, {\it Des. Codes Cryptogr.}, 86(11): 2629-2643, 2018.



\bibitem{KR2002}
D. R. Kuhn  and M. J. Reilly,  An investigation of the applicability of design of experiments to software testing, Proceedings of the 27th NASA/ IEEE Software Engineering Workshop, NASA Goddard Space Flight Center, 91-95, 2002.



\bibitem{KW2004}
D. R. Kuhn and D. R. Wallace, Software fault interaction and implication for software testing, {\it IEEE Trans. Softw. Eng.}, 30(6): 1-4, 2004.

\bibitem{RMS2018}
S. Raaphorst, L. Moura and B. Stevens, Variable strength covering arrays, {\it J. Combin. Des.}, 26(9): 417-438, 2018.

\bibitem{R1947}
C. R. Rao, Factorial experiments derivable from combinatorial arrangements of arrays, Supplement to the Journal of the Royal Statistical Society, 9: 128-139, 1947.



\bibitem{STY2012}
C. Shi, Y. Tang and J. Yin,  The equivalence between optimal detecting arrays and super-simple OAs, {\it Des. Codes Cryptogr.}, 62: 131-142, 2012.


\bibitem{SY2014}
C. Shi and J. Yin, Existence of super-simple OA$_\lambda(3,5,v)]'$s, {\it Des. Codes Cryptogr.}, 72: 369-380, 2014.

\bibitem{TY2011}
Y. Tang and J. Yin, Detecting arrays and their optimality, {\it Acta Math. Sin., Engl. Ser.}, 27: 2309-2318, 2011.


\bibitem{TMPS2017}
G. Tzanakis, L. Moura,  D. Panario and B. Stevens, Covering arrays from $m$-sequences and character sums,  {\it Des. Codes Cryptogr.}, 85(3): 437-456, 2017.



\bibitem{WNXS2007}
Z. Wang, C. Nie, B. Xu and L. Shi, Optimal Test Suite Generation Methods for Neighbor Factors Combinatorial Testing (in chinese), {\it Chinese Journal of Computers.}, 30(2): 200-211, 2007.



\end {thebibliography}
\end{document}